\documentclass[conference]{IEEEtran}
\IEEEoverridecommandlockouts
\usepackage{cite}
\usepackage{amsmath,amssymb,amsfonts,amsthm}
\usepackage{graphicx}
\usepackage{textcomp}
\usepackage{xcolor}
\def\BibTeX{{\rm B\kern-.05em{\sc i\kern-.025em b}\kern-.08em
    T\kern-.1667em\lower.7ex\hbox{E}\kern-.125emX}}

\usepackage{enumerate}
      
\usepackage{bm} 
\usepackage{subfig}

\usepackage{algorithm}
\usepackage[noend]{algpseudocode} 

\usepackage{tikz}
\usepackage{pgfplots}

\definecolor{orange}{RGB}{217,83,25}
\definecolor{blue}{RGB}{0,114,189}
\definecolor{yellow}{RGB}{237,177,32}
\definecolor{purple}{RGB}{126,47,142}
\definecolor{green}{RGB}{119,172,48}

\newcommand{\Cline}[2]{\raisebox{2pt}{\tikz{\draw[#1,#2,line width=1.5pt](0,0) -- (5mm,0);}}}
\newcommand{\Ccross}[1]{\raisebox{-1pt}{\tikz{\draw[color=#1,solid,line width=2pt](-1mm,-1mm) -- (1mm,1mm) (-1mm,1mm) -- (1mm,-1mm);}}}

\graphicspath{{Figures/}}%

\usepackage{layouts}

\newtheorem{remark}{Remark}
\newtheorem{definition}{Definition}

\newtheorem{lem}{Lemma}
\newtheorem{thm}{Theorem}
\newtheorem*{pf*}{Proof}

\usepackage{eso-pic}
\AddToShipoutPictureBG*{%
	\AtPageUpperLeft{%
		\setlength\unitlength{1in}%
		\hspace*{\dimexpr0.5\paperwidth\relax}
		\makebox(0,-1.75)[c]{
			\begin{tabular}{c c}
				S.A.N. Nouwens \emph{et al.},
				Constraint-Adaptive MPC for linear systems:\\ A system-theoretic framework for speeding up MPC through online constraint removal, \\
				Accepted for
				{\em Automatica},
				uploaded to ArXiv July 13, 2023 \\
\end{tabular}}}}

\begin{document}

\title{\vspace{2em}Constraint-Adaptive MPC for linear systems: A system-theoretic framework for speeding up MPC through online constraint removal\\
	\thanks{This research is supported by KWF Kankerbestrijding and NWO Domain AES, as part of their joint strategic research programme: Technology for Oncology II. The collaboration project is co-funded by the PPP Allowance made available by Health$\sim$Holland, Top Sector Life Sciences \& Health, to stimulate public-private partnerships.\\
	\\
	© 2023. This manuscript version is made available under the CC-BY-NC-ND 4.0 license https://creativecommons.org/licenses/by-nc-nd/4.0/}
}

\author{\IEEEauthorblockN{1\textsuperscript{st} S.A.N. Nouwens}
	\IEEEauthorblockA{\textit{Department of Mechanical Engineering} \\
		\textit{Eindhoven University of Technology}\\
		s.a.n.nouwens@tue.nl}
	\and
	\IEEEauthorblockN{2\textsuperscript{nd} M.M. Paulides}
	\IEEEauthorblockA{\textit{Department of Electrical Engineering} \\
		\textit{Eindhoven University of Technology}}
	\IEEEauthorblockA{\textit{Department of Radiotherapy} \\
		\textit{Erasmus UMC Cancer Institute}}
	\and
	\IEEEauthorblockN{3\textsuperscript{rd} W.P.M.H. Heemels}
	\IEEEauthorblockA{\textit{Department of Mechanical Engineering} \\
		\textit{Eindhoven University of Technology}}
}

\maketitle

\begin{abstract}
Reducing the computation time of model predictive control (MPC) is important, especially for systems constrained by many state constraints. In this paper, we propose a new online constraint removal framework for linear systems, for which we coin the term constraint-adaptive MPC (ca-MPC). In so-called {\em exact} ca-MPC, we adapt the imposed constraints by removing, at each time-step, a subset of the state constraints in order to reduce the computational complexity of the receding-horizon optimal control problem, while ensuring that the closed-loop behavior is {\em identical} to that of the original MPC law. We also propose an {\em approximate} ca-MPC scheme in which a further reduction of computation time can be accomplished by a tradeoff with closed-loop performance, while still preserving recursive feasibility, stability, and constraint satisfaction properties. The online constraint removal exploits fast backward and forward reachability computations combined with optimality properties.
\end{abstract}


\section{Introduction}
Model predictive control (MPC) is a successful control technology adopted in many application fields \cite{Mayne2014,Mayne2000}, and is based on recursively solving a finite-horizon optimization problem online. Solving an optimization problem at each time-step can prohibit the real-time feasibility of the controller for computationally complex scenarios. This is particularly the case in applications requiring the control of systems with many state constraints, which is the setting studied in this paper. 

Efforts to improve the computational aspects of MPC are commonplace in the literature with, amongst others, explicit MPC, model reduction, and tailored numerical solvers as prominent examples, see, e.g., \cite{Arnstrom2022May,Genuit2011,Bemporad2002,Bemporad2011,Frison2020Jan,Hovland2006,Jerez2011,Rawlings2019}. In particular, constraint removal techniques were developed to accelerate MPC for systems subject to many constraints. 

Constraint removal techniques can be roughly separated into offline \cite{Ardakani2015,Paulraj2010,Roald2019}, and online methods \cite{Jost2013,Jost2015}. Although of interest, offline methods can be prohibitively complex to compute and do not always enable real-time MPC, as constraints can only be removed if they are redundant for \emph{all} feasible states. The fact that offline methods do not depend on and thus can not exploit the current state information might make them less effective than online constraint removal techniques. 

In contrast, online techniques {\em can} exploit knowledge of the current state and possibly even more. As a result, online techniques have the potential to remove considerably more constraints compared to offline methods. For example, in \cite{Jost2013} so-called regions of activity for each constraint are (approximately) computed that are based on the initial state for the MPC problem. Loosely speaking, the region of activity represents the set of initial states for which the particular constraint in consideration is active at the minimizer of the MPC problem. The online complexity of this method is reported to scale linearly in the number of constraints. However, approximating the region of activity can become intractable in scenarios with many constraints, as for each constraint either an ellipsoidal or hypercube outer approximation must be computed that depends on all other constraints. Alternatively in \cite{Jost2015}, a Lyapunov-based approach is proposed, assuming the cost function of the MPC problem is a Lyapunov function. Here, for each inequality constraint, the cost function is minimized assuming the particular inequality is active in the sense of equality. This leads to the minimum cost function value for which the particular constraint can be active. Next, when the MPC control law is running, the value of the Lyapunov function for the current state is compared to the pre-computed values for all constraints. If this value for the current state is lower than the corresponding value for an inequality constraint, it can be removed permanently from the MPC problem. Interestingly, both the method based on the regions of activity and the method using a Lyapunov cost function are complementary to the framework we will present in this paper.

In this work, we will present a new online constraint removal framework for linear systems, called {\em constraint-adaptive} MPC (ca-MPC). We present both {\em exact} and {\em approximate} ca-MPC strategies. Crucially, in exact ca-MPC, the closed-loop behavior of the resulting accelerated MPC feedback law is {\em identical} to that of the original MPC feedback law. In approximate ca-MPC a further reduction of computation time can be accomplished compared to exact ca-MPC due to a tradeoff with closed-loop performance (as the closed-loop behavior is no longer identical in approximate ca-MPC). However, in approximate ca-MPC crucial properties such as recursive feasibility, stability, and constraint satisfaction can be preserved by design. Both strategies exploit system-theoretic properties, such as reachability and optimality, in a computationally effective manner. The method presented in this paper extends our preliminary work in \cite{Nouwens2021} and \cite{Nouwens2021b} in which only initial ideas were presented (without any technical proofs). The current work formalizes these initial ideas in a complete framework and specifies also the technical underlying results and their rigorous proofs. A new numerical case study, extending the earlier one, is provided as well, which shows a two-order reduction in computational time of the ca-MPC scheme compared to the original MPC scheme, while still having identical closed-loop behavior. We also show how the preliminary results of \cite{Nouwens2021} can be seen as a special case, {\em approximate} ca-MPC (see Section~\ref{sec:journal_further_extension}).


\section{System and MPC setup}\label{sec:journal_preliminaries}
In this paper, we consider plants that can be described by a discrete-time linear time-invariant (LTI) system
\begin{align}\label{eq:journal_system_description}
	&\bm{x}_{k+1} = \bm{Ax}_k + \bm{Bu}_k,
\end{align}
although several ideas also apply to nonlinear and time-varying plants, see, e.g., \cite{Nouwens2021}. In \eqref{eq:journal_system_description}, $\bm{x}_k\in\mathbb{R}^n$ and $\bm{u}_k\in\mathbb{R}^m$ denote the plant states and the inputs, respectively, at discrete time $k\in\mathbb{N}$. Furthermore, $\bm{A}\in\mathbb{R}^{n\times n}$ and $\bm{B}\in\mathbb{R}^{n\times m}$. The system \eqref{eq:journal_system_description} is subject to polyhedral state and input constraints given for $k\in\mathbb{N}$ by
\begin{subequations}\label{eq:journal_state_and_input_constraints}
	\begin{align} 
		\bm{x}_k &\in \mathbb{X} := \{\bm{x}\in\mathbb{R}^n\mid \bm{c}_j\bm{x}\leq b_j,\ \text{for}\ j\in\mathbb{N}_{[1,n_x]}\},\\
		\bm{u}_k &\in \mathbb{U} := \{\bm{u}\in\mathbb{R}^m\mid \bm{g}_j\bm{u} \leq h_j,\ \text{for}\ j\in\mathbb{N}_{[1,n_u]}\}. \label{eq:journal_input_constraints}
	\end{align}
\end{subequations}
Here, $\mathbb{X}$ and $\mathbb{U}$ are assumed to be non-empty polyhedral sets with $\bm{c}_j\in\mathbb{R}^{1\times n}$, $\bm{g}_j\in\mathbb{R}^{1\times m}$, $b_j\in\mathbb{R}$, and $h_j\in\mathbb{R}$. In this paper, we study systems that are constrained by many state constraints, i.e., $n_x\gg1$.

\subsection{MPC setup} 
Based on the system dynamics \eqref{eq:journal_system_description} and constraints \eqref{eq:journal_state_and_input_constraints}, a common MPC setup, given state $\bm{x}_k$ at time $k\in\mathbb{N}$, is
\begin{subequations}\label{eq:journal_basic_mpc}
	\begin{align}
		\underset{\bm{X}_k,\ \bm{U}_k}{\text{minimize}}\ \quad&J(\bm{X}_k,\bm{U}_k),\label{eq:journal_basic_mpc_a}\\
		\text{subject to }\quad&\bm{X}_k = \bm{\Phi x}_k+\bm{\Gamma U}_k ,\label{eq:journal_basic_mpc_b}\\
		&\bm{X}_k \in \mathcal{X}:=\textstyle\prod_{i=1}^{N}\mathbb{X}_i,\label{eq:journal_basic_mpc_c}\\
		&\bm{U}_k \in \mathcal{U}:=\mathbb{U}^{N},\label{eq:journal_basic_mpc_d}
	\end{align}
	where
	\begin{align}
		&J(\bm{X}_k,\bm{U}_k) := \ell_T(\bm{x}_{N|k}) + \textstyle\sum_{i=0}^{N-1}\ell(\bm{x}_{i|k},\bm{u}_{i|k}),\label{eq:journal_basic_mpc_costfun}\\
		&\bm{X}_k:=[\bm{x}_{1|k}^\top \cdots \bm{x}_{N|k}^\top]^\top, \bm{U}_k:=[\bm{u}_{0|k}^\top \cdots \bm{u}_{N-1|k}^\top]^\top,\label{eq:journal_basic_mpc_g}\\
		&\mathbb{X}_i:=\{\bm{x}\in\mathbb{R}^n\mid\bm{c}_{i,j}\bm{x}\leq b_{i,j},\ \text{for}\ j\in\mathbb{N}_{[1,n_{x_i}]}\},\label{eq:journal_basic_mpc_h}\\
		&\bm{\Phi} =	\left[\begin{smallmatrix}
			\bm{A} \\
			\bm{A}^2 \\
			\vdots \\
			\bm{A}^N
		\end{smallmatrix}\right],\qquad \bm{\Gamma}:=\left[\begin{smallmatrix}
			\bm{B} & \bm{0} & \cdots & \bm{0}\\ 
			\bm{AB} & \bm{B} & \cdots & \bm{0} \\
			\vdots & \vdots & \ddots & \vdots \\
			\bm{A}^{N-1}\bm{B} & \bm{A}^{N-2}\bm{B} &\cdots &\bm{B} 
		\end{smallmatrix}\right].
	\end{align}
\end{subequations}
Here, $\ell$, $\ell_T$, $\bm{x}_{i|k}$, $\bm{u}_{i|k}$, and $\mathbb{X}_i$ denote the stage cost, the terminal cost, the predicted state, the predicted input, and the state constraint set at predicted time $i\in\mathbb{N}_{[1,N]}:=\{1,2,\cdots,N\}$ made at time $k\in\mathbb{N}$, respectively. The state constraints depend on $i$ for generality and to facilitate compact notation using the Cartesian product. Typically, $\mathbb{X}_i$ is chosen as $\mathbb{X}_i=\mathbb{X}\subset\mathbb{R}^n$ for $i\in\mathbb{N}_{[1,N-1]}$ and $\mathbb{X}_N=\mathbb{X}_T\subseteq\mathbb{X}$, where $\mathbb{X}_T$ denotes a suitable controlled invariant terminal set \cite{Mayne2000}. The $i|k$ subscript is used to denote the $i$-th prediction at time $k$.

For the optimization problem \eqref{eq:journal_basic_mpc}, we denote the set of feasible input sequences parameterized by $\bm{x}_k$ as
\begin{align}\label{eq:journal_U_fn}
	\mathcal{U}_{f}(\bm{x}_k) :=\{ \bm{U}_k \in \mathcal{U} \mid \eqref{eq:journal_basic_mpc_b},\ \eqref{eq:journal_basic_mpc_c} \},
\end{align}
and the set of feasible states by $\mathbb{X}_f :=\{ \bm{x} \in \mathbb{X} \mid \mathcal{U}_{f}(\bm{x}) \neq \emptyset \}$. Under suitable assumptions on $\ell$, $\ell_T$, $\mathcal{U}$, and $\mathcal{X}$, e.g., $\ell$ and $\ell_T$ being continuous and $\mathcal{X}$ being closed and $\mathcal{U}$ being compact \cite{Mayne2000}, a minimizer of \eqref{eq:journal_basic_mpc} exists for all $\bm{x}_k\in \mathbb{X}_f$ and we denote by $\bm{U}^\star_{k} := [\bm{u}^{\star\top}_{0|k}\ \cdots\ \bm{u}^{\star\top}_{N-1|k}]^\top$ a particular one at time $k\in\mathbb{N}$ for state $\bm{x}_k$, i.e., 
\begin{align}\label{eq:journal_argmin}
	\bm{U}^\star_{k} \in \mathcal{U}^\star(\bm{x}_k) := \underset{\bm{U}_k \in \mathcal{U}_{f}(\bm{x}_k)}{\arg \min} \bar{J}(\bm{x}_k,\bm{U}_k),
\end{align}
where $\bar{J}(\bm{x}_k,\bm{U}_k) := J(\bm{\Phi x}_k+\bm{\Gamma U}_k,\bm{U}_k)$. The set of all optimal predicted state sequences corresponding to $\mathcal{U}^\star(\bm{x}_k)$ is denoted by
\begin{align}
	\mathcal{X}^\star(\bm{x}_k) := \bm{\Phi x}_k + \bm{\Gamma} \mathcal{U}^\star(\bm{x}_k),
\end{align}
where a particular one is given by $\bm{X}_k^\star=\bm{\Phi x}_k + \bm{\Gamma} \bm{U}^\star_{k}\in\mathcal{X}^\star(\bm{x}_k)$. Using a receding horizon implementation, the MPC problem \eqref{eq:journal_basic_mpc} is turned into a feedback law $K_\text{MPC}:\mathbb{X}_f\rightarrow\mathbb{U}$ by applying the first computed input in $\bm{U}_k^\star$ on the real plant \eqref{eq:journal_system_description}, i.e., $\bm{u}_k := K_\text{MPC}(\bm{x}_k):= \bm{u}^\star_{0|k}$, $k\in\mathbb{N}$.

\subsection{Reduced MPC problem}
To address the problem of removing redundant state constraints from \eqref{eq:journal_basic_mpc}, we introduce the reduced MPC problem, where the original constraints set $\mathcal{X}$ is replaced by a (state-dependent) reduced constraint set denoted by $\mathcal{X}^\text{red}(\mathcal{A}(\bm{x}_k))=\textstyle\prod_{i=1}^N\mathbb{X}^\text{red}_i(\mathbb{A}_i(\bm{x}_k))$, leading to
\begin{subequations}\label{eq:journal_red_mpc}
	\begin{align}
		\underset{\bm{X}_k,\ \bm{U}_k}{\text{minimize}}\ \quad&J(\bm{X}_k,\bm{U}_k),\label{eq:journal_red_mpc_a}\\
		\text{subject to }\quad &\eqref{eq:journal_basic_mpc_b},\ \eqref{eq:journal_basic_mpc_d},\label{eq:journal_red_mpc_b}\\
		&\bm{X}_k \in \mathcal{X}^\text{red}(\mathcal{A}(\bm{x}_k)).\label{eq:journal_red_mpc_c}
	\end{align}
\end{subequations}
The reduced constraint set $\mathcal{X}^\text{red}(\mathcal{A}(\bm{x}_k))$ is described by an index-set $\mathcal{A}(\bm{x}_k)=\prod_{i=1}^N\mathbb{A}_i(\bm{x}_k)$, such that
\begin{align}\label{eq:journal_def_red}
	\mathbb{X}^\text{red}_i(\mathbb{A}_i) &:= \{\bm{x}\in\mathbb{R}^n\mid \bm{c}_{i,j}\bm{x}\leq b_{i,j},\ \text{for }j\in\mathbb{A}_i\}.
\end{align}
Clearly, $\mathcal{X}^\text{red}$ is defined by a subset of the constraints in $\mathcal{X}$. Moreover, observe that the index set $\mathcal{A}(\bm{x}_k)$ depends on $\bm{x}_k$. Hence, our ca-MPC scheme requires the specification of the set-valued mapping $\mathbb{A}_i:\mathbb{X}_f\rightrightarrows\mathbb{N}_{[1,n_{x_i}]}$ for $i\in\mathbb{N}_{[1,N]}$. A formal problem formulation will be given in Section~\ref{sec:journal_exact_ca-MPC}. We use the notation $\rightrightarrows$ to indicate the set-valuedness of the maps of $\mathbb{A}_i$ in the sense that $\mathbb{A}_i(\bm{x}_k)\subseteq\mathbb{N}_{[1,n_{x_i}]}$ for $i\in\mathbb{N}_{[1,N]}$.

The reduced MPC problem \eqref{eq:journal_red_mpc} with the reduced constraint sets \eqref{eq:journal_def_red} gives rise to the set of minimizers
\begin{subequations}
	\begin{align}\label{eq:journal_argmin_red}
		&\mathcal{U}^{\text{red}\star}(\bm{x}_k,\mathcal{A}(\bm{x}_k)) := \underset{\bm{U}_k \in \mathcal{U}_{f}^\text{red}(\bm{x}_k,\mathcal{A}(\bm{x}_k))}{\arg \min} \hspace{-1em}\bar{J}(\bm{x}_k,\bm{U}_k), \\ \label{eq:journal_U_fn_red}	
		&\mathcal{U}_{f}^\text{red}(\bm{x}_k,\mathcal{A}(\bm{x}_k)) :=\{ \bm{U}_k \in \mathcal{U} \mid \eqref{eq:journal_red_mpc_b}-\eqref{eq:journal_red_mpc_c}\},\\
		&\mathbb{X}_f^\text{red}:=\{\bm{x}\in\mathbb{X}\mid\mathcal{U}_{f}^\text{red}(\bm{x}_k,\mathcal{A}(\bm{x}_k))\neq\emptyset\}.
	\end{align}
\end{subequations}
Again, similar to the original MPC problem, we introduce the set of ``optimal'' state trajectories
\begin{align}
	\mathcal{X}^{\text{red}\star}(\bm{x}_k,\mathcal{A}(\bm{x}_k)):=\bm{\Phi x}_k + \bm{\Gamma}\mathcal{U}^{\text{red}\star}(\bm{x}_k,\mathcal{A}(\bm{x}_k)).
\end{align}

\subsection{Preliminaries and notation}
Given a set $\mathbb{V}\subset\mathbb{R}^n$, we denote the affine transformation of $\mathbb{V}$ with matrix $\bm{M}\in\mathbb{R}^{m\times n}$ and vector $\bm{b}\in\mathbb{R}^m$ by $	\bm{M}\mathbb{V}+\bm{b} := \{\bm{Mv}+\bm{b}\in\mathbb{R}^{m}\mid \bm{v}\in\mathbb{V}\}$.
We define an ellipsoidal set $\mathcal{E}(\bm{L},\bm{q}):=\{\bm{x}\in\mathbb{R}^n\mid \|\bm{L}(\bm{x}-\bm{q})\|_2\leq 1\}$, where $\bm{L}\in\mathbb{R}^{n\times n}$, $\bm{q}\in\mathbb{R}^n$, with $\|\bm{v}\|_2^2:=\sum_{i=1}^nv_i^2$ for $\bm{v}\in\mathbb{R}^n$. To project from sets defined for state (or input) trajectories to individual state vectors (or input vectors), we introduce the projection $\text{P}_i$ of a set $\mathcal{V}\subseteq\mathbb{R}^{Nn}$ as
\begin{align}
	\text{P}_i(\mathcal{V}) := \{\bm{v}_i\in\mathbb{R}^n\mid &\exists \bm{v}_j\in\mathbb{R}^n,\ j\in\mathbb{N}_{[1,N]}\backslash\{i\}\\ \nonumber &\text{s.t. } [\bm{v}_1^\top\ \cdots\ \bm{v}_N^\top]^\top\in\mathcal{V}\},
\end{align}
where $i\in\mathbb{N}_{[1,N]}$, assuming $N$ is clear from the context. We define the normal cone of a set $\mathcal{V}\subset\mathbb{R}^d$ at $\bm{v}\in\mathcal{V}$ as $N_\mathcal{V}(\bm{v}):=\{\bm{p}\in\mathbb{R}^{d}\mid \bm{p}^\top(\bm{y}-\bm{v})\leq 0\ \text{for all}\ \bm{y}\in\mathcal{V}\}$.


\section{Exact constraint-adaptive MPC}\label{sec:journal_exact_ca-MPC}
We now introduce \emph{exact ca-MPC}, in which the term ``exact" refers to the property that the closed-loop system will \emph{not} be changed when replacing \eqref{eq:journal_basic_mpc} by the reduced MPC problem \eqref{eq:journal_red_mpc}. 
\begin{definition}\label{def:journal_exact_campc}
	A ca-MPC scheme based on \eqref{eq:journal_red_mpc} for given reduced constraint mappings $\mathbb{A}_i:\mathbb{X}_f\rightrightarrows\mathbb{N}_{[1,n_{x_i}]}$, $i\in\mathbb{N}_{[1,N]}$, is called exact, if $\mathbb{X}_f=\mathbb{X}_f^\text{red}$ and the set of minimizers of the reduced MPC problem \eqref{eq:journal_red_mpc} is the same as for the original MPC problem \eqref{eq:journal_basic_mpc} for all $\bm{x}_k\in\mathbb{X}_f$, i.e.,
	\begin{align}\label{eq:journal_campc_exactness}
		\mathcal{U}^{\star}(\bm{x}_k) = 	\mathcal{U}^{\text{red}\star}(\bm{x}_k,\mathcal{A}(\bm{x}_k)) .
	\end{align}
\end{definition}
As the minimizers are identical, \emph{exact} ca-MPC trivially inherits all performance, stability, and constraint satisfaction properties from the original MPC problem.

\emph{The main problem} considered in this paper can now be formulated as follows: Construct computationally tractable set-valued mappings $\mathbb{A}_i:\mathbb{X}_f\rightrightarrows\mathbb{N}_{[1,n_{x_i}]}$, with the number of constraints $\text{card}(\mathbb{A}_i(\bm{x}_k))\ll n_{x_i}$ (if possible) for $i\in\mathbb{N}_{[1,N]}$ for all $\bm{x}_k\in\mathbb{X}_f$, such that, the resulting ca-MPC scheme is exact in the sense of Definition~\ref{def:journal_exact_campc}.

By creating simpler MPC problems, ca-MPC can accelerate both interior-point and active-set solvers. For interior-point methods, ca-MPC straightforwardly reduces the complexity of each Newton step, thereby accelerating the optimization problem. Additionally, ca-MPC is also a natural extension to active-set solvers, as the working and inactive constraint sets are already dynamically updated. By removing constraints a priori, determining which constraint has to be added to the working set is easier, as there are fewer constraints to evaluate.

\emph{The key concept} of exact ca-MPC is inspired by the following basic observation. A ca-MPC scheme based on \eqref{eq:journal_red_mpc} for given mappings $\mathcal{A}$, is exact if and only if 	
\begin{align}\label{eq:journal_key_concept_prop}
	\mathcal{X}^{\text{red}\star}(\bm{x}_k,\mathcal{A}(\bm{x}_k))\subseteq \mathcal{X},
\end{align}
for all $\bm{x}_k\in\mathbb{X}_f$. Indeed, when the optimal set of state trajectories satisfy all state constraints, we obtain the same set of minimizers. However, note that \eqref{eq:journal_key_concept_prop} is not suitable as a direct tool to design appropriate set-valued mappings $\mathcal{A}$, as \eqref{eq:journal_key_concept_prop} cannot be used in a constructive manner. To overcome this, we build upon the following theorem using outer approximations of $\mathcal{X}^{\text{red}\star}(\bm{x}_k,\mathcal{A}(\bm{x}_k))$ by a set $\mathcal{M}(\bm{x}_k)$ (not depending on $\mathcal{A}$).
\begin{thm}\label{thm:journal_thm1}
	Consider the original and reduced MPC optimization problem \eqref{eq:journal_basic_mpc} and \eqref{eq:journal_red_mpc}, respectively, with $\mathbb{A}_i:\mathbb{X}_f\rightrightarrows\mathbb{N}_{[1,n_{x_i}]}$ and $\mathcal{M}:\mathbb{X}_f\rightrightarrows\mathbb{R}^{Nn}$. If for all $\bm{x}_k\in\mathbb{X}_f$,
	\begin{align}\tag{C1}\label{eq:journal_thm1_cond1}
		\mathcal{X}^{\text{red}\star}(\bm{x}_k,{\mathcal{A}}(\bm{x}_k))&\subseteq\mathcal{M}(\bm{x}_k),\\ \tag{C2}\label{eq:journal_thm1_cond2}
		\mathcal{M}(\bm{x}_k) \cap \textstyle\mathcal{X}^\text{red}(\mathcal{A}(\bm{x}_k)) &\subseteq \mathcal{X},
	\end{align}
	then $\mathbb{X}_f=\mathbb{X}_f^\text{red}$ and the ca-MPC scheme \eqref{eq:journal_red_mpc} is exact.
\end{thm}
\begin{pf*}
	Let $\bm{x}_k\in\mathbb{X}_f$ be given.	First, observe
	\begin{align}\label{eq:journal_thm1_pf1}
		\mathcal{U}_{f}(\bm{x}_k) \subseteq \mathcal{U}^\text{red}_{f}(\bm{x}_k,\mathcal{A}(\bm{x}_k)),
	\end{align}
	which is immediate as \eqref{eq:journal_red_mpc} uses a subset of the constraints of \eqref{eq:journal_basic_mpc}. Hence, $\mathbb{X}_f\subseteq\mathbb{X}_f^\text{red}$ and for all $\bm{x}_k\in\mathbb{X}_f$
	\begin{align}\label{eq:journal_thm1_pf2}
		\underset{\bm{U}_k \in \mathcal{U}_{f}(\bm{x}_k)}{\min} \bar{J}(\bm{x}_k,\bm{U}_k) \geq \underset{\bm{U}_k \in \mathcal{U}^\text{red}_{f}(\bm{x}_k,\mathcal{A}(\bm{x}_k))}{\min} \hspace{-1em}\bar{J}(\bm{x}_k,\bm{U}_k).
	\end{align} 
	Next, to show \eqref{eq:journal_campc_exactness}, take $\bm{x}_k\in\mathbb{X}_f^\text{red}$, for which each
	\begin{align}\label{eq:journal_thm1_pf3}
		\bm{U}^{\text{red}\star}_k \in \mathcal{U}^{\text{red}\star}(\bm{x}_k,\mathcal{A}(\bm{x}_k))
	\end{align} 
	satisfies $\bm{U}^{\text{red}\star}_k\in\mathcal{U}$ by \eqref{eq:journal_red_mpc_c} and $\bm{X}_{k}^{\text{red}\star}\in\mathcal{X}$ by combining \eqref{eq:journal_red_mpc_b}, \eqref{eq:journal_thm1_cond1}, and \eqref{eq:journal_thm1_cond2}. Hence, we obtain $\mathcal{U}^{\text{red}\star}(\bm{x}_k,\mathcal{A}(\bm{x}_k))\in\mathcal{U}_f(\bm{x}_k)$, and by extension, $\mathbb{X}_f^\text{red}\subseteq\mathbb{X}_f$. Therefore, we have, $\mathbb{X}_f^\text{red}=\mathbb{X}_f$ and
	\begin{align}\label{eq:journal_thm1_pf4}
		\underset{\bm{U}_k \in \mathcal{U}_{f}(\bm{x}_k)}{\min} \bar{J}(\bm{x}_k,\bm{U}_k) = \underset{\bm{U}_k \in \mathcal{U}^\text{red}_{f}(\bm{x}_k,\mathcal{A}(\bm{x}_k))}{\min} \hspace{-1em}\bar{J}(\bm{x}_k,\bm{U}_k).
	\end{align} 
	Combining \eqref{eq:journal_thm1_pf2}, \eqref{eq:journal_thm1_pf3}, and \eqref{eq:journal_thm1_pf4} leads now to \eqref{eq:journal_campc_exactness}.\qed
\end{pf*}
The main takeaway from Theorem~\ref{thm:journal_thm1} is that we concentrate the $\mathcal{A}$-dependence in $\mathcal{X}^\text{red}(\mathcal{A}(\bm{x}_k))$. The way we apply Theorem~\ref{thm:journal_thm1} is to first construct an $\mathcal{A}$-independent set $\mathcal{M}(\bm{x}_k)$ satisfying \eqref{eq:journal_thm1_cond1}, for all possible choices of $\mathcal{A}$. Once, this $\mathcal{M}(\bm{x}_k)$ is available, we can select, in the second step, the indices in $\mathcal{A}(\bm{x}_k)$ such that \eqref{eq:journal_thm1_cond2} is satisfied. Intuitively, when $\mathcal{M}$ satisfies \eqref{eq:journal_thm1_cond1}, then, its elements indicate which state constraints can be violated. Therefore, by adding these state constraint indices to $\mathcal{A}$, we ensure that these ``at-risk" constraints are not violated. As a result, the optimal state trajectory $\mathcal{X}^{\text{red}\star}(\bm{x}_k,{\mathcal{A}}(\bm{x}_k))$ will satisfy all state constraints, and, hence, we obtain an exact ca-MPC scheme.

An important observation is that a smaller $\mathcal{M}$ has the potential to remove more constraints, see \eqref{eq:journal_thm1_cond2}. However, $\mathcal{M}$ still has to satisfy \eqref{eq:journal_thm1_cond1}. To this end, a useful extension to Theorem~\ref{thm:journal_thm1}, is to utilize additional information besides the initial state, see, e.g., Sections \ref{sec:journal_backward_reach_set} and \ref{sec:journal_cost_level_set}. For example, we can exclude certain state constraints from the removal process in order to find a smaller $\mathcal{M}$. To this end, we introduce a set of \emph{fixed constraints} given by $\mathcal{F}(\bm{x}_k) := \prod_{i=1}^N\mathbb{F}_i(\bm{x}_k)$, where the index set $\mathcal{A}$ must satisfy $\mathcal{F}(\bm{x}_k)\subseteq\mathcal{A}(\bm{x}_k)$. We will omit the dependence of $\mathcal{M}$ on $\mathcal{F}$, as it will be clear from context.

Besides introducing fixed constraints to exclude constraints from the removal process, it is useful to construct $\mathcal{A}$ in a subtractive manner. By doing so, we can be conservative in removing constraints without losing exactness, but this will provide many computational benefits. To make this concrete, we parameterize $\mathcal{A}$ as
\begin{align}
	\mathbb{A}_i(\bm{x}_k) = \mathbb{N}_{[1,n_{x_i}]}\backslash\mathbb{I}_i(\bm{x}_k),\quad i\in\mathbb{N}_{[1,N]},
\end{align}
where $\mathbb{I}_i:\mathbb{X}_f\rightrightarrows\mathbb{N}_{[1,n_{x_i}]}$, with $\mathbb{I}(\bm{x}_k)\cap\mathbb{F}_i(\bm{x}_k)=\emptyset$ for all $\bm{x}_k\in\mathbb{X}_f$, denotes the set of removed constraints from the MPC problem. We also define the compact notation $\mathcal{I}(\bm{x}_k):=\prod_{i=1}^N\mathbb{I}_i(\bm{x}_k)$. 
\begin{lem}\label{lem:journal_I}
	Consider the original and reduced MPC problems \eqref{eq:journal_basic_mpc} and \eqref{eq:journal_red_mpc}, respectively, and let $\mathbb{F}_{i}:\mathbb{X}_f\rightrightarrows\mathbb{N}_{[1,n_{x_i}]}$ and $\mathcal{M}:\mathbb{X}_f\rightrightarrows\mathbb{R}^{Nn}$ satisfying \eqref{eq:journal_thm1_cond1} be given. Let the set-valued mapping $\mathbb{I}_i:\mathbb{X}_f\rightrightarrows\mathbb{N}_{[1,n_{x_i}]}$, with $\mathbb{I}_i(\bm{x}_k)\cap\mathbb{F}_i(\bm{x}_k)=\emptyset$, $i\in\mathbb{N}_{[1,N]}$, capture the removed state constraints, i.e., $\mathbb{A}_i(\bm{x}_k) = \mathbb{N}_{[1,n_{x_i}]}\backslash\mathbb{I}_i(\bm{x}_k), i\in\mathbb{N}_{[1,N]}$. If for all $\bm{x}_k\in\mathbb{X}_f$
	\begin{align}\label{eq:journal_lem_F_I}\tag{C3}
		\mathcal{M}(\bm{x}_k) \cap \mathcal{X}^\text{red}(\mathcal{I}(\bm{x}_k)) = 	\mathcal{M}(\bm{x}_k).
	\end{align}
	Then, the resulting ca-MPC scheme is exact.
\end{lem}
\begin{pf*}
	To show \eqref{eq:journal_thm1_cond2}, we start with the inclusion $\mathcal{X}\cap\mathcal{M}(\bm{x}_k)\subseteq\mathcal{X}$. Since $\mathbb{A}_i\cup\mathbb{I}_i = \mathbb{N}_{[1,n_{x_i}]}$, $i\in\mathbb{N}_{[1,N]}$, we obtain $\mathcal{X}^\text{red}(\mathcal{A}(\bm{x}_k))\cap\mathcal{X}^\text{red}(\mathcal{I}(\bm{x}_k)) \cap \mathcal{M}(\bm{x}_k) \subseteq\mathcal{X}$. Exploiting \eqref{eq:journal_lem_F_I} gives $\mathcal{X}^\text{red}(\mathcal{A}(\bm{x}_k))\cap \mathcal{M}(\bm{x}_k) \subseteq\mathcal{X}$, which is equivalent to \eqref{eq:journal_thm1_cond2}. \qed
\end{pf*}
The result of Lemma~\ref{lem:journal_I} is instrumental as it provides a computationally efficient method to remove constraints from the reduced MPC problem. For example, \eqref{eq:journal_lem_F_I} can be evaluated for each constraint independently. As a result, we can also consider projections of $\mathcal{M}$, i.e., 
\begin{align}
	\text{P}_i(\mathcal{M}(\bm{x}_k)) \cap \mathbb{X}^\text{red}_i(\mathbb{I}_i(\bm{x}_k))=\text{P}_i(\mathcal{M}(\bm{x}_k)),
\end{align}
for all $i\in\mathbb{N}_{[1,N]}$, is equivalent to \eqref{eq:journal_lem_F_I}. This observation will lead to an efficient ca-MPC implementation, as we will see in the next section.


\section{Proposed exact ca-MPC implementation}\label{sec:journal_proposed_implementation}
In this section, we will make the exact ca-MPC scheme concrete by introducing three sets $\mathcal{M}^{(1)}$, $\mathcal{M}^{(2)}$, and $\mathcal{M}^{(3)}$ to construct $\mathcal{M}$ as $\mathcal{M}^{(1)}\cap\mathcal{M}^{(2)}\cap\mathcal{M}^{(3)}$. Hereafter, we will provide a concrete outline on how $\mathcal{M}$ is used to compute the mapping $\mathcal{A}$.

\subsection{Forward reachable set}\label{sec:journal_forward_reach_set}
As all state trajectories start at $\bm{x}_k$ and satisfy the system dynamics and input constraints, we can use the input-constrained forward reachable set to construct $\mathcal{M}^{(1)}(\bm{x}_k)$. To this end, we introduce the \emph{input-constrained} forward reachable set $\overrightarrow{\mathcal{H}}(\bm{x}_k,\mathbb{U}):=\prod_{i=1}^N\overrightarrow{\mathbb{H}}_i(\bm{x}_k,\mathbb{U})\subset\mathbb{R}^{Nn}$, which can be computed by the recursion
\begin{subequations}\label{eq:journal_fwd}
	\begin{align}
		\overrightarrow{\mathbb{H}}_{i+1}(\bm{x}_k,\mathbb{U}) :=& \bm{A}\overrightarrow{\mathbb{H}}_{i}(\bm{x}_k,\mathbb{U}) + \bm{B}\mathbb{U},\\
		\overrightarrow{\mathbb{H}}_{0}(\bm{x}_k,\mathbb{U}) =& \{\bm{x}_k\}.
	\end{align}
\end{subequations} 
To compute the state-dependent forward reachable set in a real-time setting we define $\mathcal{M}^{(1)}(\bm{x}_k)$ as
\begin{align}\label{eq:journal_M_fwd}
	\mathcal{M}^{(1)}(\bm{x}_k)=\textstyle\prod_{i=1}^N\mathcal{E}(\bm{L}_{i,1},\bm{q}_{i,1}) + \bm{A}^i\bm{x}_k,
\end{align} 
where the ellipsoidal set satisfies $\mathcal{E}(\bm{L}_{i,1},\bm{q}_{i,1})\supset\overrightarrow{\mathbb{H}}_{i}(\bm{0},\mathbb{U})$, which can be computed using \cite{Halder2018,henk2012lowner}. The construction \eqref{eq:journal_M_fwd} solves two problems. First, note that the forward reachable set can be decomposed as $\overrightarrow{\mathbb{H}}_{i}(\bm{x}_k,\mathbb{U}) = \bm{A}^{i}\bm{x}_k + \overrightarrow{\mathbb{H}}_{i}(\bm{0},\mathbb{U})$ for linear systems. Hence, we can compute the forward reachable set once offline for a zero initial state, and shift the result by the free state response $\bm{A}^i\bm{x}_k$. Second, we use ellipsoids that have a fixed complexity to outer approximate the true forward reachable set, which can be arbitrarily complex.

\subsection{Backward reachable set}\label{sec:journal_backward_reach_set}
Given the forward reachable set, a natural extension is to exploit backward reachability too. Indeed, when we choose $\mathcal{F}(\bm{x}_k) = \prod_{i=1}^{N-1}\emptyset\times\mathbb{N}_{[1,n_{x_N}]}$, all optimal state trajectories that can result from \eqref{eq:journal_red_mpc} must end in the terminal set. To this end, we introduce the \emph{input- and terminal-state-constrained} backward reachable set for \eqref{eq:journal_system_description} and \eqref{eq:journal_input_constraints}, i.e., $\overleftarrow{\mathcal{H}}(\mathbb{X}_N,\mathbb{U}):=\prod_{i=1}^N\overleftarrow{\mathbb{H}}_{i}(\mathbb{X}_N,\mathbb{U})\subset\mathbb{R}^{Nn}$, which can be computed using the recursion
\begin{subequations}\label{eq:journal_bwd}
	\begin{align}\nonumber 
		\overleftarrow{\mathbb{H}}_{i-1}(\mathbb{X}_N,\mathbb{U}) := &\{\bm{x}\in\mathbb{R}^n\mid\bm{Ax}+\bm{Bu} \in \overleftarrow{\mathbb{H}}_{i}(\mathbb{X}_N,\mathbb{U})\\
		&\text{for some}\ \bm{u}\in\mathbb{U}\},\ i\in\mathbb{N}_{[2,N]}\\
		\overleftarrow{\mathbb{H}}_{N}(\mathbb{X}_N,\mathbb{U}) = &\mathbb{X}_N.
	\end{align}
\end{subequations}
Note that the backwards reachable set is not state-dependent but can be prohibitively complex, similar to the forward reachable set. Hereto, we define $\mathcal{M}^{(2)}$ as
\begin{align}
	\mathcal{M}^{(2)}(\bm{x}_k)=\textstyle\prod_{i=1}^N\mathcal{E}(\bm{L}_{i,2},\bm{q}_{i,2}), 
\end{align}
where $\mathcal{E}(\bm{L}_{i,2},\bm{q}_{i,2})\subseteq\overleftarrow{\mathbb{H}}_{i}(\mathbb{X}_N,\mathbb{U})$ to manage online complexity.

\subsection{First-order optimality set}\label{sec:journal_cost_level_set}
The third set, $\mathcal{M}^{(3)}$, exploits the optimization based nature of MPC. For example, a feasible solution to \eqref{eq:journal_basic_mpc} upper bounds the cost function, thereby bounding the minimizer \cite{Nouwens2021b}. One method to obtain a feasible solution, is by extending the minimizer $\bm{U}^\star_{k-1}$ obtained at time $k-1$ (as is common when using terminal set and cost methods, see \cite{Mayne2000}). In this section, we introduce a more advanced notion of this concept by exploiting the property that optimal state trajectories satisfy first-order optimality conditions, which, in combination with a feasible solution to \eqref{eq:journal_basic_mpc}, yields a smaller set. 

\begin{figure}[th!]
	\centering
	\includegraphics[width=8.2cm]{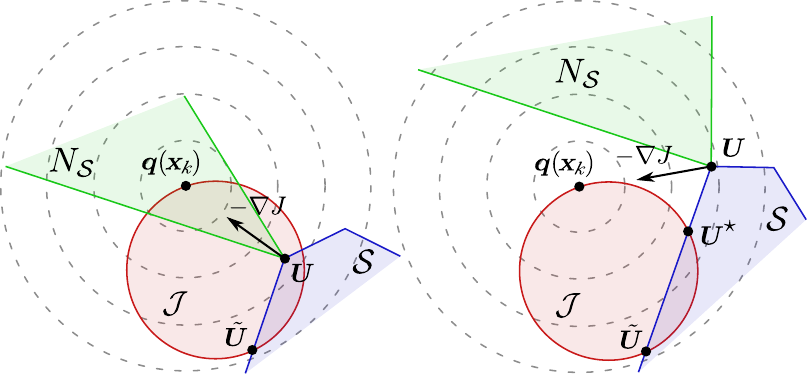}
	\caption{Two illustrations of \eqref{eq:journal_cost_levelset} for a quadratic cost function at different points $\bm{U}$. The dashed lines denote level sets of the cost function and $\bm{q}(\bm{x}_k)$ denotes the unconstrained minimizer, where $\nabla_{\bm{U}}\bar{J}(\bm{x}_k,\bm{q}(\bm{x}_k)) = \bm{0}$. Note that $\mathcal{J}$ is significantly smaller than the level set of $\bar{J}$ at $\tilde{\bm{U}}$. \emph{Left}: $\bm{U}\in\mathcal{J}$, as there exists a set $\mathcal{S}$ such that $N_\mathcal{S}$ includes $-\nabla_{\bm{U}}\bar{J}$. In fact, $\bm{U}$ is the minimizer. \emph{Right}: $\bm{U}\notin\mathcal{J}$, as there does not exist a set $\mathcal{S}$ such that $-\nabla_{\bm{U}}\bar{J}\in N_\mathcal{S}$. As expected, we observe that the minimizer for $\mathcal{S}$ satisfies $\bm{U}^\star\in\mathcal{J}$.}
	\label{fig:journal_costlevelset}
\end{figure}

We start by formalizing the first-order optimality conditions that optimal input trajectories satisfy, under differentiability of $\bar{J}$ with respect to $\bm{U}$,
\begin{align}\label{eq:journal_first_order_opt}
	-\nabla_{\bm{U}} \bar{J}(\bm{x}_k,\bm{U}_k)\in N_{\mathcal{U}^\text{red}_{f}(\bm{x}_k,\mathcal{A}(\bm{x}_k))}(\bm{U}^\star),
\end{align}
where $N_{\mathcal{U}^\text{red}_{f}(\bm{x}_k,\mathcal{A}(\bm{x}_k))}(\bm{U}^\star)$ denotes the normal cone of $\mathcal{U}^\text{red}_{f}(\bm{x}_k,\mathcal{A}(\bm{x}_k))$ at $\bm{U}^\star$, as is commonly seen in stationary conditions for optimization problems \cite[Sec.~12.7]{Nocedal2006}. Observe that \eqref{eq:journal_first_order_opt} is not a constructive condition to design $\mathcal{M}^{(3)}$, as it uses $\mathcal{A}$ and $\bm{U}^\star$ in its definition. To get around this issue, we introduce the set of input trajectories for which a convex constraint set exists, such that, first-order optimality conditions can be satisfied given a feasible input trajectory $\tilde{\bm{U}}_k\in\mathcal{U}_f(\bm{x}_k)$, see Figure~\ref{fig:journal_costlevelset}. Note that the set of all convex sets that include the feasible input trajectory is guaranteed to include $\mathcal{U}^\text{red}_{f}(\bm{x}_k,\mathcal{A}(\bm{x}_k))$. Hereto, we define $\mathcal{J}(\bm{x}_k,\tilde{\bm{U}}_k)\subset\mathbb{R}^{Nm}$, 
\begin{align}\label{eq:journal_cost_levelset}
	&\mathcal{J}(\bm{x}_k,\tilde{\bm{U}}_k):=\{\bm{U}\in\mathbb{R}^{Nm}\mid \exists\text{convex}\ \mathcal{S}\subseteq\mathbb{R}^{Nm},\\ \nonumber
	&\text{s.t. }\bm{U},\tilde{\bm{U}}_k\in\mathcal{S},\ \text{and }-\nabla_{\bm{U}}\bar{J}(\bm{x}_k,\bm{U}_k)\in N_\mathcal{S}(\bm{U})\}.
\end{align} 
Given \eqref{eq:journal_cost_levelset}, we obtain $\mathcal{M}^{(3)}(\bm{x}_k) = \bm{\Phi}\bm{x}_k + \bm{\Gamma}\mathcal{J}(\bm{x}_k,\tilde{\bm{U}}_k)$. 

While \eqref{eq:journal_cost_levelset} seems complex at first sight, it has an elegant solution for quadratic cost functions. To show this, we specify a quadratic cost function $J$ as in \eqref{eq:journal_basic_mpc_costfun} with $\ell(\bm{x},\bm{u})=\bm{x}^\top\bm{Q}\bm{x}+\bm{u}^\top\bm{R}\bm{u},\ \ell_T(\bm{x})=\bm{x}^\top\bm{P}\bm{x}$, where $\bm{Q}$ and $\bm{P}$ denote positive semi-definite matrices and $\bm{R}$ is positive definite. When substituting the quadratic cost and system dynamics into \eqref{eq:journal_basic_mpc_costfun}, it is easy to see that an equivalent formulation for the cost is
\begin{align}\label{eq:journal_quadratic_cost_2}
	\bar{J}(\bm{x}_k,\bm{U}_k) = &\|\bm{G}(\bm{U}_k-\bm{q}(\bm{x}_k))\|_2^2 + r(\bm{x}_k)
\end{align}
with $\bm{G}\in\mathbb{R}^{Nm\times Nm}$ and linear maps $\bm{q}:\mathbb{X}_f\rightarrow\mathbb{R}^{Nm}$, $r:\mathbb{X}_f\rightarrow\mathbb{R}$ that depend on the system dynamics, the current state, and $\bm{Q},\ \bm{P},\ \bm{R}$. Observe that $\bm{q}(\bm{x}_k)$ denotes the unconstrained minimizer of the MPC problem, which is the unique solution to $\nabla_{\bm{U}}\bar{J}(\bm{x}_k,\bm{q}(\bm{x}_k)) = \bm{0}$. Using \eqref{eq:journal_quadratic_cost_2}, $\mathcal{J}(\bm{x}_k,\tilde{\bm{U}}_k)$ is given by the ellipsoid
\begin{align}\label{eq:journal_cost_level_set_implementation}
	&\mathcal{J}(\bm{x}_k,\tilde{\bm{U}}_k)=\{\bm{U}\in\mathbb{R}^{Nm}\mid\\\nonumber
	&\|\bm{G}\big(\bm{U} - \frac{1}{2}(\tilde{\bm{U}}_k + \bm{q}(\bm{x}_k))\big)\|_2\leq\frac{1}{2}\|\bm{G}(\tilde{\bm{U}}_k - \bm{q}(\bm{x}_k))\|_2\}.
\end{align}
The proof of this result is provided in Appendix~\ref{apx:journal_cost_levelset}, and a schematic illustration is in Figure~\ref{fig:journal_costlevelset}. 

\subsection{Integration into an exact ca-MPC scheme}\label{sec:journal_integration}
Now we will present how the sets $\mathcal{M}^{(i)}$, $i=\{1,2,3\}$, are used to construct an exact ca-MPC scheme. Recall that $\mathcal{M}^{(2)}$ requires $\mathcal{F}(\bm{x}_k) = \prod_{i=1}^{N-1}\emptyset \times \mathbb{N}_{[1,n_{x_N}]}$ and $\mathcal{M}^{(3)}$ requires availability of an $\tilde{\bm{U}}_k\in\mathcal{U}_f(\bm{x}_k)$. As $\mathcal{M}$ is defined as the intersection of ellipsoids, we utilize \eqref{eq:journal_lem_F_I} to construct $\mathcal{I}(\bm{x}_k)$ as
\begin{subequations}\label{eq:journal_I_set_def}
	\begin{align}
		\mathbb{I}_i(\bm{x}_k) &= \mathbb{I}^{(1)}_i(\bm{x}_k)\cup\mathbb{I}^{(2)}_i(\bm{x}_k)\cup\mathbb{I}^{(3)}_i(\bm{x}_k),\\ \label{eq:journal_I_set_def_b}
		\mathbb{I}^{(l)}_i(\bm{x}_k) &= \{j\in\mathbb{N}_{[1,n_{x_i}]}\backslash\mathbb{F}_i(\bm{x}_k)\mid\\ \nonumber &\|\bm{c}_{i,j}\bm{L}_{i,l}^{-1}\|_2\leq|b_{i,j}-\bm{c}_{i,j}\bm{q}_{i,l}|\},
	\end{align}
\end{subequations}
where $\mathcal{E}(\bm{L}_{i,l},\bm{q}_{i,l}) = \text{P}_i(\mathcal{M}^{(l)}(\bm{x}_k))$. Observe that for each $l\in\{1,2,3\}$, it holds that 
\begin{align}\label{eq:journal_I_set_part}
	\text{P}_i(\mathcal{M}^{(l)}(\bm{x}_k)) \cap \mathbb{X}_i(\mathbb{I}^{(l)}_i(\bm{x}_k))= \text{P}_i(\mathcal{M}^{(l)}(\bm{x}_k)).
\end{align}
Hence, \eqref{eq:journal_I_set_def} implies \eqref{eq:journal_lem_F_I}. The derivation of \eqref{eq:journal_I_set_def_b} is given in Appendix~\ref{apx:journal_set_intersection}. Crucially, $\bm{L}_{i,l}$ is not state-dependent (up to a constant), therefore, $\|\bm{c}_{i,j}\bm{L}_{i,l}^{-1}\|_2$ can be pre-computed (up to a constant) for all $i,j,l$. As a result, the ellipsoidal description of $\mathcal{M}$ allows for the computationally efficient of $\mathcal{I}(\bm{x}_k)$ using \eqref{eq:journal_I_set_def}. In fact, \eqref{eq:journal_I_set_def_b} only requires one inner product, a scalar absolute value, and a comparison.

For the remainder of this section, we present an overview of the integrated ca-MPC scheme in Algorithm~\ref{alg:journal_caMPC_concrete}. In Algorithm~\ref{alg:journal_caMPC_concrete}, we use $\leftarrow$ to denote ``compute using". First, in lines~\ref{algl:journal_caMPC1}-\ref{algl:journal_caMPC5}, in an offline setting, we pre-compute the forward and backward reachable sets and corresponding ellipsoidal outer approximations. Second, in an online setting, we start by computing the first-order optimality set in line~\ref{algl:journal_caMPC7} based on the measured state $\bm{x}_k$ and generated feasible input sequence $\tilde{\bm{U}}_k$. Next, in lines~\ref{algl:journal_caMPC_start}-\ref{algl:journal_caMPC_end}, we perform the constraint removal using the forward and backward reachable sets and first-order optimality set. Finally, we compute the minimizer of the resulting reduced MPC problem (line \ref{algl:journal_caMPC_computeU}) and apply $\bm{u}_k = \bm{u}_{0|k}^{\text{red}\star}$ to the plant (line \ref{algl:journal_caMPC_applyU}). 
\begin{algorithm}
	\caption{Implementation of exact ca-MPC}
	\begin{algorithmic}[1]
		\State $k=0,\ N\in\mathbb{N}$\label{algl:journal_caMPC1}
		\State $\overrightarrow{\mathbb{H}}_i(\bm{0},\mathbb{U})\leftarrow$ \eqref{eq:journal_fwd}, $i\in\mathbb{N}_{[1,N]}$
		\State $\overleftarrow{\mathbb{H}}_i(\mathbb{X}_N,\mathbb{U})\leftarrow$ \eqref{eq:journal_bwd}, $i\in\mathbb{N}_{[1,N-1]}$
		\State $\bm{L}_{i,1},\ \bm{q}_{i,1} \leftarrow \overrightarrow{\mathbb{H}}_i(\bm{0},\mathbb{U})$, s.t. \eqref{eq:journal_fwd}, $i\in\mathbb{N}_{[1,N]}$
		\State $\bm{L}_{i,2},\ \bm{q}_{i,2} \leftarrow \overleftarrow{\mathbb{H}}_i(\mathbb{X}_N,\mathbb{U})$, s.t. \eqref{eq:journal_bwd}, $i\in\mathbb{N}_{[1,N-1]}$\label{algl:journal_caMPC5}
		\While{\textsc{true}}
		\State \textsc{Measure }$\bm{x}_k\in\mathbb{X}_f$
		\State $\mathcal{J}(\bm{x}_k,\tilde{\bm{U}}_k)\leftarrow$ \eqref{eq:journal_cost_levelset} \textsc{with} $\tilde{\bm{U}}_k\in\mathcal{U}_f(\bm{x}_k)$\label{algl:journal_caMPC7}
		
		\For{$i=1,2,\cdots,N-1$}	\label{algl:journal_caMPC_start}
		\State $\mathbb{I}_i^{(1)} \leftarrow$ \eqref{eq:journal_I_set_def} \textsc{with} $\mathcal{E}(\bm{L}_{i,1}, \bm{q}_{i,1})+\bm{A}^i\bm{x}_k$	\label{algl:journal_caMPC12}
		\State $\mathbb{I}_i^{(2)} \leftarrow$ \eqref{eq:journal_I_set_def} \textsc{with} $\mathcal{E}(\bm{L}_{i,2}, \bm{q}_{i,2})$	
		\State $\mathbb{I}_i^{(3)} \leftarrow$ \eqref{eq:journal_I_set_def} \textsc{with} $\text{P}_i(\bm{\Gamma}\mathcal{J}(\bm{x}_k,\tilde{\bm{U}}_k)) + \bm{A}^i\bm{x}_k$			
		\State $\mathbb{A}_i\leftarrow\mathbb{N}_{[1,n_{x_i}]}\backslash(\mathbb{I}_i^{(1)}\cup\mathbb{I}_i^{(2)}\cup\mathbb{I}_i^{(3)})$	\label{algl:journal_caMPC15}
		\EndFor		\label{algl:journal_caMPC_end}	
		
		\State $\bm{U}_k^{\text{red}\star}\leftarrow$ \eqref{eq:journal_argmin_red} \label{algl:journal_caMPC_computeU}
		\State \textsc{Apply control input} $\bm{u}_k\leftarrow\bm{u}_{0|k}^{\text{red}\star}$ \label{algl:journal_caMPC_applyU}
		\State $k\leftarrow k+1$
		\EndWhile		
	\end{algorithmic}\label{alg:journal_caMPC_concrete}
\end{algorithm}


\section{Approximate ca-MPC}\label{sec:journal_further_extension}
The exact ca-MPC method introduced in Sections~\ref{sec:journal_exact_ca-MPC} and \ref{sec:journal_proposed_implementation} is highly effective in the sense that the online complexity of computing $\mathcal{A}(\bm{x}_k)$ is low, while significant reductions in the number of constraints can be achieved. A relaxation, \emph{approximate ca-MPC}, aims to remove even more constraints from the MPC problem. Approximate ca-MPC will maintain the constraint satisfaction (and stability), but will no longer be exact.

Based on \eqref{eq:journal_thm1_cond2}, we observe that shrinking $\mathcal{M}$ can result in the removal of more state constraints. One way to realize this, is by designing a smaller $\mathcal{M}$ that satisfies \eqref{eq:journal_thm1_cond1}. Hereto, we impose \emph{additional} input constraints that allow us to design a tighter $\mathcal{M}$. To this end, we augment the original MPC problem \eqref{eq:journal_basic_mpc} with the input constraint	$\bm{U}_k\in\tilde{\bm{U}}_k + \delta\mathbb{U}^N$, where $\tilde{\bm{U}}_k\in\mathcal{U}_f(\bm{x}_k)$ is a feasible input sequence, and $\delta\mathbb{U}\subset\mathbb{R}^m$ with $\bm{0}\in\delta\mathbb{U}$ is a set limiting the input sequence $\bm{U}_k$ to be ``close" to $\tilde{\bm{U}}_k$. By constraining the minimizer to be close to a feasible input trajectory, we obtain ``smaller" forward and backward reachable sets that can lead to the removal of more constraints from the MPC problem. We will consider a particular choice of $\tilde{\bm{U}}_k$, namely, $\tilde{\bm{U}}_k = [\bm{u}_{1|k-1}^{\star\top}\ \cdots\ \bm{u}_{N-1|k-1}^{\star\top}\ \tilde{\bm{u}}^{\top}(\bm{x}_{N|k-1}^\star)]^\top$, 
where $\tilde{\bm{u}}:\mathbb{X}_N\rightarrow\mathbb{U}$ denotes an auxiliary feedback law that renders $\mathbb{X}_N$ positively invariant, i.e., $\bm{Ax}+\bm{B}\tilde{\bm{u}}(\bm{x})\in\mathbb{X}_N$ for all $\bm{x}\in\mathbb{X}_N$. The MPC problem corresponding to this new setup is given by
\begin{subequations}\label{eq:journal_basic_mpc_extended}
	\begin{align}
		\underset{\bm{X}_k,\ \bm{U}_k}{\text{minimize}}\ \quad&J(\bm{X}_k,\bm{U}_k),\label{eq:journal_basic_mpc_extendedc_a}\\
		\text{subject to }\quad&\eqref{eq:journal_basic_mpc_b},\ \eqref{eq:journal_basic_mpc_d} \label{eq:journal_basic_mpc_extended_b}\\
		&\bm{U}_k \in \mathcal{U}\cap(\tilde{\bm{U}}_k + \delta\mathbb{U}^N).\label{eq:journal_basic_mpc_extended_c}
	\end{align}
\end{subequations}
Note that depending on $\tilde{\bm{u}}:\mathbb{X}_N\rightarrow\mathbb{U}$, the closed-loop stability of \eqref{eq:journal_basic_mpc_extended} can be unchanged with respect to \eqref{eq:journal_basic_mpc}, as the standard stability proof based on a terminal set and cost \cite{Mayne2000}, is still applicable. Interestingly, the approximate ca-MPC scheme is exact with respect to \eqref{eq:journal_basic_mpc_extended}, but not \eqref{eq:journal_basic_mpc}, hence, its properties can be analyzed using \eqref{eq:journal_basic_mpc_extended}. Last, note that the size of $\delta\mathbb{U}$ is a tuning knob (and can even depend on time $k$ and prediction step $i$) that can trade-off potential performance loss to increased constraint removal, and thus computational benefits.


\section{Numerical example}\label{sec:journal_numerical_example}
\begin{figure}[th!]
	\sbox0{\Cline{black}{solid}}\sbox1{\Cline{black}{dashed}}\sbox2{\Ccross{orange}}
	\centering
	\includegraphics[width=6cm]{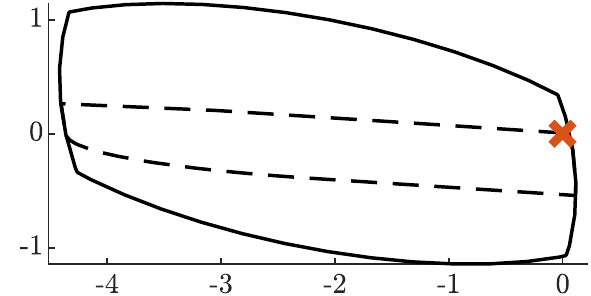}
	\caption{Illustration of the state constraints $\mathbb{X}_i$, $i\in\mathbb{N}_{[1,N-1]}$ (\usebox0), the terminal set $\mathbb{X}_N$ (\usebox1), and the origin (\usebox2).}
	\label{fig:journal_constraint_sets}
\end{figure}

\begin{figure}[th!]
	\sbox0{\Cline{purple}{solid}}\sbox1{\Cline{yellow}{dashed}}\sbox2{\Cline{blue}{dashed}}\sbox3{\Ccross{orange}}
	\centering
	\includegraphics[width=5cm]{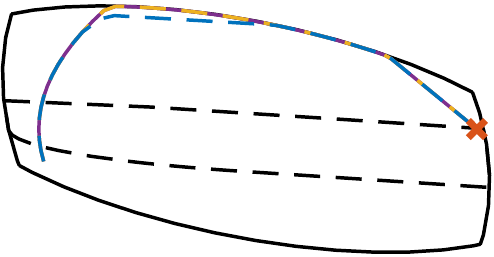}
	\caption{Closed-loop trajectories for the original MPC (\usebox0), exact ca-MPC (\usebox1), and approximate ca-MPC (\usebox2).}
	\label{fig:journal_trajectory}
\end{figure}

\begin{figure}[th!]
	\sbox0{\Cline{black}{solid}}\sbox1{\Cline{black}{dashed}}
	\centering
	\includegraphics[width=8.4cm]{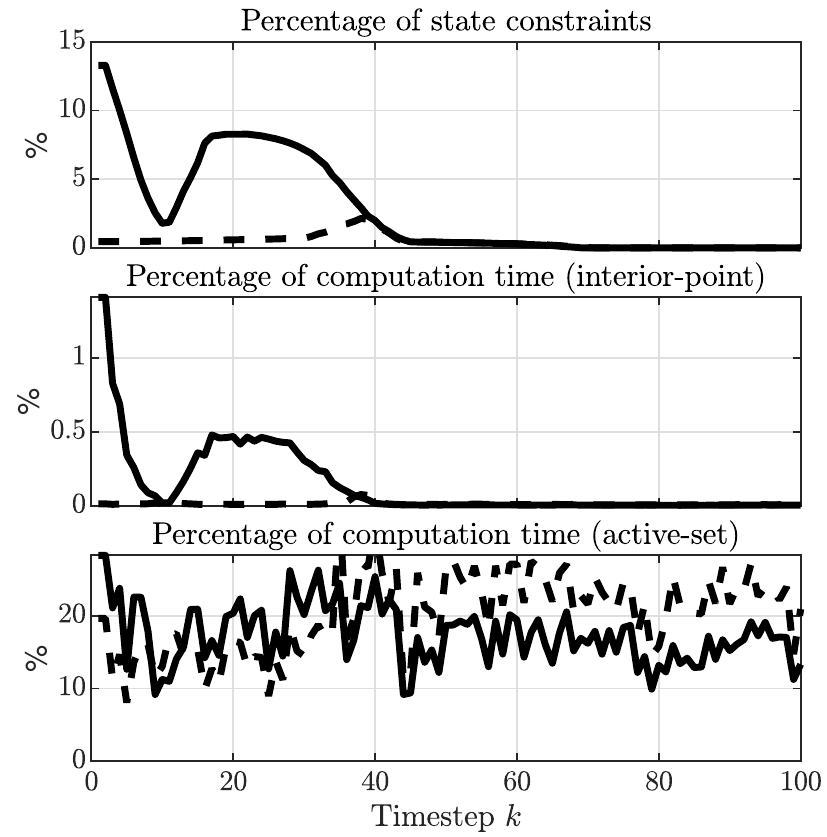}
	\caption{The percentage of state constraints and computation time for exact ca-MPC (\usebox0) and approximate ca-MPC (\usebox1) relative to the original MPC feedback law for both an interior-point and active-set solver.}
	\label{fig:journal_percentage}
\end{figure}

We will use a double integrator system
\begin{align}\label{eq:journal_example_system}
	\bm{x}_{k+1} = \left[\begin{smallmatrix}
		1 & 0.1 \\ 0 & 1
	\end{smallmatrix}\right]\bm{x}_k + \left[\begin{smallmatrix}
		0.005 \\ 0.1
	\end{smallmatrix}\right]u_k,
\end{align}
as the two-dimensional state allows for convenient visualization. The input constraints are $\mathbb{U}=\{u\in\mathbb{R}\mid |u|\leq1\}$. To obtain an example with a large number of linear state constraints as in \eqref{eq:journal_basic_mpc}, we approximate two quadratic constraints with linear inequalities, i.e.,
\begin{align}\label{eq:journal_double_integrator}
	\hspace{-0.2em}\mathbb{X}_i:=&\{\bm{x}\mid (\bm{v}_{1,j}-\bm{d})^\top\bm{P}_1(\bm{x}-\bm{d})\leq 1,\ j\in\mathbb{N}_{[1,n_v]}\}\\ \nonumber
	\cap&\{\bm{x}\mid (\bm{v}_{2,j}-\bm{d})^\top\bm{P}_2(\bm{x}-\bm{d})\leq 1,\ j\in\mathbb{N}_{[1,n_v]}\},
\end{align}
for $i\in\mathbb{N}_{[1,N-1]}$, $\bm{d}=[2.15\ 0]^\top$, and 
\begin{subequations}
	\begin{align}
		&\bm{P}_1 = \left[\begin{smallmatrix}
			0.14 &  0.17\\
			0.17 & 0.97
		\end{smallmatrix}\right],\quad \bm{P}_2 = \left[\begin{smallmatrix}
			0.20 & 0.05\\
			0.05 & 0.21
		\end{smallmatrix}\right], \\
		&(\bm{v}_{1,j}-\bm{d})^\top\bm{P}_1(\bm{v}_{1,j}-\bm{d})=1,\\ &(\bm{v}_{2,j}-\bm{d})^\top\bm{P}_2(\bm{v}_{2,j}-\bm{d})=1,\ j\in\mathbb{N}_{[1,n_v]}.
	\end{align}
\end{subequations}
The points $\bm{v}_{1,j},\ \bm{v}_{2,j},\ j\in\mathbb{N}_{[1,n_v]}$, are chosen on the boundary of the respective ellipses, which will result in a tangent inequality constraint at $\bm{v}_{1,j}$ and $\bm{v}_{2,j}$, respectively. Note that all constraints are non-redundant, i.e., the removal of each constraint changes the set. Moreover, the number of state constraints will be $n_{x_i} = 2 n_v$, $i\in\mathbb{N}_{[1,N-1]}$ (for the results in Figures~\ref{fig:journal_trajectory} and \ref{fig:journal_percentage} we picked $n_v = 330$, which results in $660$ state constraints for $i\in\mathbb{N}_{[1,N-1]}$). The terminal set $\mathbb{X}_N\subseteq\mathbb{X}_1$ is chosen to be positively invariant for the control law $u_k = \bm{K}_T\bm{x}_k=-[0.01\ 0.01]\bm{x}_k$. An illustration of the state constraints, including the terminal set, is shown in Figure~\ref{fig:journal_constraint_sets}. When choosing $N=12$, the total number of state constraints in this example is $\sum_{i=1}^{N}n_{x_i} = 2(N-1)n_v + n_{x_N} = 7468$.

We define the MPC cost function as $\ell(\bm{x},\bm{u})=\bm{x}^\top\bm{Q}\bm{x}+\bm{u}^\top\bm{R}\bm{u},\ \ell_T(\bm{x})=\bm{x}^\top\bm{P}\bm{x}$ using $\bm{Q} = \left[\begin{smallmatrix}
	1 & 0 \\ 0 & 1
\end{smallmatrix}\right],\ \bm{P} = \left[\begin{smallmatrix}
	1 & 0 \\ 0 & 1
\end{smallmatrix}\right],\ \bm{R} = 1
$. Both the first-order optimality set and the approximate ca-MPC scheme use the feasible input sequence, $\tilde{\bm{U}}_k = [\bm{u}_{1|k-1}^{\star\top} \cdots \bm{u}_{N-1|k-1}^{\star\top}\ (\bm{K}_T\bm{x}^\star_{N|k-1})^\top
]^\top$. The last component needed for the \emph{approximate} ca-MPC scheme is the additional input constraint, which we take as
$\delta\mathbb{U} := \{u\in\mathbb{R}\mid |u|\leq 0.3\}$.

\begin{figure}[th!]
	\sbox0{\Ccross{blue}}\sbox1{\Cline{yellow}{solid}}\sbox2{\Cline{purple}{solid}}\sbox3{\Cline{green}{solid}}
	\centering
	\includegraphics[width=5cm]{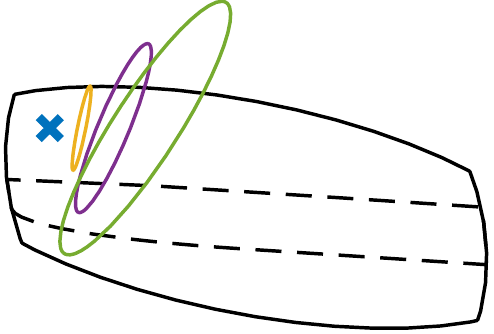}
	\caption{Illustration of the outer ellipsoidal approximations of the forward reachable sets $\mathcal{E}(\bm{L}_{i,1},\bm{q}_{i,1})+\bm{A}^i\bm{x}_k$ for $i=4$ ({\usebox1}), $i=8$ ({\usebox2}), $i=12$ ({\usebox3}), starting from $\bm{x}_k$ ({\usebox0}).}
	\label{fig:journal_fwd_sets}
\end{figure}

\begin{figure}[th!]
	\sbox0{\Cline{black}{dashed}}\sbox1{\Cline{yellow}{solid}}\sbox2{\Cline{purple}{solid}}\sbox3{\Cline{green}{solid}}
	\centering
	\includegraphics[width=5cm]{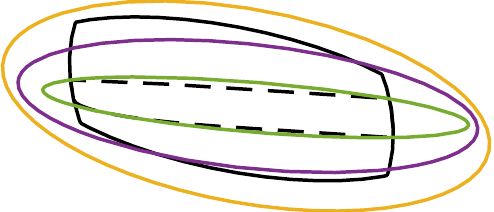}
	\caption{Illustration of the outer ellipsoidal approximations of the backward reachable sets $\mathcal{E}(\bm{L}_{i,2},\bm{q}_{i,2})$ for $i=4$ ({\usebox1}), $i=8$ ({\usebox2}), $i=12$ ({\usebox3}), starting from the $\mathbb{X}_N$ (\usebox0).}
	\label{fig:journal_bwd_sets}
\end{figure}

\begin{figure}[!ht]
	\sbox0{\Ccross{blue}}\sbox1{\Cline{yellow}{solid}}\sbox2{\Cline{purple}{solid}}\sbox3{\Cline{green}{solid}}
	\centering
	\subfloat[]{
		\includegraphics[width=4cm]{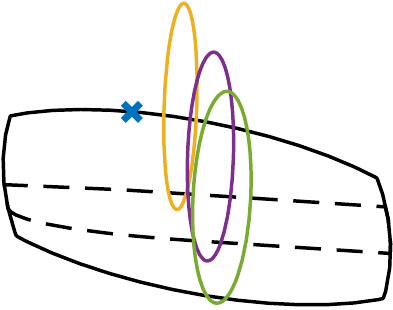}
	}
	\subfloat[]{
		\includegraphics[width=4cm]{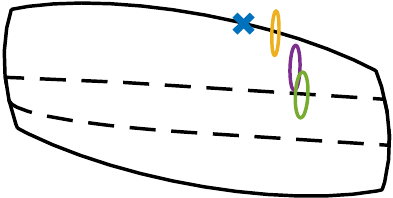}
	}
	\caption{Illustration of the first-order optimality sets $\text{P}_i(\bm{\Gamma}\mathcal{J}(\bm{x}_k,\tilde{\bm{U}}_k)) + \bm{A}^i\bm{x}_k$ for $i=4$ ({\usebox1}), $i=8$ ({\usebox2}), $i=12$ ({\usebox3}) for two initial states $\bm{x}_k$ (\usebox0). The ellipsoids shrink when $\bm{x}_k$ approaches the origin, i.e., $\tilde{\bm{U}}_k$ approaches $\bm{q}(\bm{x}_k)$.}
	\label{fig:journal_j_sets}
\end{figure}

The two-dimensional state of \eqref{eq:journal_example_system} allows for the visualization of $\mathcal{M}^{(l)}(\bm{x}_k)$ for $l\in\{1,2,3\}$, see Figures~\ref{fig:journal_fwd_sets}-\ref{fig:journal_j_sets}. Note that forward and backward reachable set are approximated using the outer Löwner-John ellipsoid of \eqref{eq:journal_fwd} and \eqref{eq:journal_bwd}, respectively. These offline ca-MPC computations required approximately 30 seconds for this example.


\begin{figure}[th!]
	\sbox0{\Cline{black}{solid}}\sbox1{\Cline{orange}{solid}}\sbox2{\Cline{blue}{solid}}\sbox3{\Cline{blue}{dashed}}
	\centering
	\includegraphics[width=8.4cm]{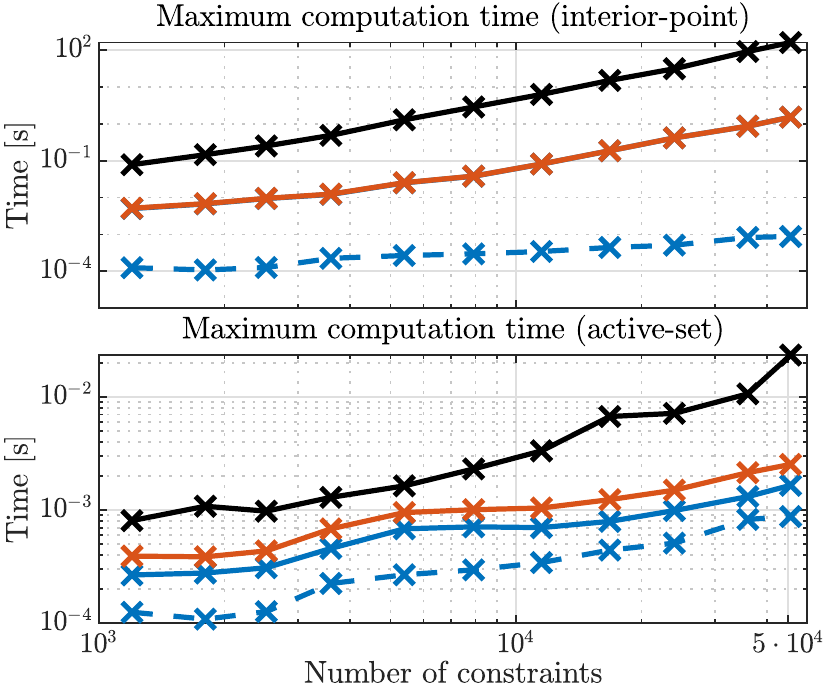}
	\caption{Maximum computation time of the original MPC (\usebox0) and exact ca-MPC scheme (\usebox1) for each solver with an increasing number of state constraints. The computation time of exact ca-MPC is split in computing $\mathcal{A}$ (\usebox3) and the resulting QP (\usebox2). Note that the QP time for the interior-point solver is virtually equal to the the total time.}
	\label{fig:journal_logTime}
\end{figure}

\emph{Simulation results:} we initialize \eqref{eq:journal_example_system} at $\bm{x}_0 = -[4\ 0.4]^\top$ and let the MPC, exact ca-MPC, and the approximate ca-MPC scheme regulate the state to the origin. The quadratic programs are solved using both a primal-dual interior-point method from the \textsc{Matlab} MPC toolbox and the DAQP active-set solver \cite{Arnstrom2022May}. The resulting state trajectories are shown in Figure~\ref{fig:journal_trajectory}. First of all, none of the schemes violated the state constraints, as expected. Second, as also expected, the exact ca-MPC scheme has an indistinguishable closed-loop trajectory compared to the original MPC solution, while, the approximate ca-MPC scheme did result in a different closed-loop trajectory. In Figure~\ref{fig:journal_percentage}, the computation time and the number of reduced constraints over time for both ca-MPC schemes and solvers are shown. Both exact and approximate ca-MPC schemes, when solved using an interior-point solver, are approximately 100-1000 times faster to compute compared to the original MPC problem. In addition, a smaller, but significant computation time improvement is observed for active-set solvers. Last, we breakdown the exact ca-MPC computation time and compare it to the original MPC problem for a range of state constraints by adjusting $n_v$ in \eqref{eq:journal_double_integrator} (total constraints vary between $1200$ and $45600$). In Figure~\ref{fig:journal_logTime}, a consistent two-order of magnitude improvement in computational time is observed for exact ca-MPC when using an interior-point solver. A smaller, but consistent, two to ten times improvement in computational time is observed when using an active-set solver. 

\begin{remark}
	While the double-integrator is convenient for illustrating all aspects of the ca-MPC scheme due to its 2-dimensional visualization possibilities, the interested reader is referred for a higher-dimensional example to \cite{Nouwens2021b}, due to space limitations. In \cite{Nouwens2021b} a similar ca-MPC scheme is applied to a discretized thermal PDE with a temperature upper bound on the (discrete) spatial domain. This leads to a discrete-time LTI model with 2000 states. Here, a similar two-order of magnitude improvement in computational time was observed for an MPC setup with 20,000 constraints.
\end{remark}


\section{Conclusions}\label{sec:journal_conclusion}
In this paper, we presented an efficient online constraint removal framework for accelerating MPC for linear systems using system-theoretic insights. A crucial aspect of our proposed method is that the closed-loop behavior and thus properties of the reduced MPC feedback law, such as stability, performance, and constraint satisfaction, remain unchanged when compared to the original MPC feedback law. We achieve this by exploiting computationally efficient bounds on the optimal state trajectory that indicate which constraints can be removed from the MPC problem. In particular, we showed that the forward and backward reachable sets and a first-order optimality set are computationally efficient and powerful tools to remove state constraints from the MPC problem. Additionally, we presented an extension, called approximate ca-MPC, that is able to trade-off closed-loop performance with the computational complexity of the resulting reduced MPC problem, while still maintaining constraint satisfaction and stability. 

The results from a numerical example show that the resulting constraint removal scheme can achieve computational speed ups of two-orders of magnitude, without loss of closed-loop performance. An alternative example with a similar ca-MPC using a thermal system described by a PDE is presented in \cite{Nouwens2021b}, where a comparable two-orders of magnitude improvement was observed as well. Moreover, many of the conceptual ideas that we exploited to get to our constraint-adaptive MPC framework can be extended to nonlinear and time-varying systems as well, see also \cite{Nouwens2021} for first steps.

\bibliographystyle{plain}    
\bibliography{literature.bib}      

\appendix

\subsection{Proof of equation \eqref{eq:journal_cost_level_set_implementation}}\label{apx:journal_cost_levelset}
\begin{lem}\label{lem:journal_lem_apx}
	Given a point $\tilde{\bm{w}}\in\mathbb{W}$, where $\mathbb{W}\subseteq\mathbb{R}^n$ is compact and convex, then $\bm{w}^\star:=\underset{\bm{w}\in\mathbb{W}}{\arg\min}\quad\|\bm{L}(\bm{w}-\bm{q})\|_2^2$ satisfies $\bm{w}^\star\in\{\bm{w}\mid-(\bm{w}-\bm{q})^\top\bm{L}^\top\bm{L}(\tilde{\bm{w}}-\bm{w})\leq 0\}$.
\end{lem}
\begin{pf*}
	First-order optimality requires $-2\bm{L}^\top\bm{L}(\bm{w}^\star-\bm{q})\in N_{\mathbb{W}}(\bm{w}^\star)$ \cite[Sec 12.7]{Nocedal2006}. Since $\tilde{\bm{w}}, \bm{w}^\star\in\mathbb{W}$, the normal cone is bounded by $N_{\mathbb{W}}(\bm{w}^\star)\subseteq\{\bm{p}\mid\bm{p}(\tilde{\bm{w}} - \bm{w}^\star)\leq0\}$. Hence, the minimizer must satisfy $-2\bm{L}^\top\bm{L}(\bm{w}^\star-\bm{q})\in\{\bm{p}\mid\bm{p}(\tilde{\bm{w}} - \bm{w}^\star)\leq0\}$. Equivalently, we obtain $-(\bm{w}^\star-\bm{q})^\top\bm{L}^\top\bm{L}(\tilde{\bm{w}} - \bm{w}^\star)\leq0$.\qed
\end{pf*}
Rewriting the result using the center of the ellipse $\frac{1}{2}(\tilde{\bm{w}} + \bm{q})$, gives $\bm{w}^\star \in \{\bm{w}\mid\|\bm{L}(\bm{w}-\frac{1}{2}(\tilde{\bm{w}} + \bm{q}))\|_2\leq\frac{1}{2}\|\bm{L}(\tilde{\bm{w}}-\bm{q})\|_2\}$, which is the same result as used in \eqref{eq:journal_cost_level_set_implementation}.

\subsection{Equivalence of \eqref{eq:journal_I_set_def_b} and \eqref{eq:journal_I_set_part}}\label{apx:journal_set_intersection}
When considering individual constraints \eqref{eq:journal_I_set_part} becomes 
\begin{align}\label{eq:journal_apx_equality}
	\mathcal{E}(\bm{L},\bm{q})\cap\{\bm{x}\in\mathbb{R}^n\mid \bm{cx}\leq b\} = \mathcal{E}(\bm{L},\bm{q}),
\end{align}
Note that $\mathcal{E}(\bm{L},\bm{q})\cap\{\bm{x}\in\mathbb{R}^n\mid \bm{cx}\leq b\} \neq \emptyset$, due to feasibility of the MPC. Hence, \eqref{eq:journal_apx_equality} is implied if the intersection between the ellipsoid and hyperplane $\{\bm{x}\in\mathbb{R}^n\mid \bm{cx}= b\}$ is either empty or a single element.

We map the ellipsoid to a unit ball on the origin using $\bm{x} = \bm{q}+\bm{L}^{-1}\bm{w}$. Substituting the mapping into the hyperplane $\bm{cx} = b$ yields $\bm{cL}^{-1}\bm{w} = b - \bm{cq}$. The distance from the hyperplane to the origin is then $\frac{|b-\bm{cq}|}{\|\bm{cL}^{-1}\|_2}$. The hyperplane intersects the ellipse when $\frac{|b-\bm{cq}|}{\|\bm{cL}^{-1}\|_2} \leq 1$. Hence, $\|\bm{cL}^{-1}\|_2 \leq |b-\bm{cq}|$ implies \eqref{eq:journal_apx_equality}, which is \eqref{eq:journal_I_set_def_b}.

\end{document}